\documentclass[letterpaper,10 pt, conference]{IEEEtran}
\IEEEoverridecommandlockouts     
\usepackage[OT1]{fontenc}
\pdfoutput=1
\usepackage[english]{babel}
\usepackage{textcomp}
\usepackage{amstext}
\usepackage{mathtools}
\usepackage{amsthm}
\usepackage{amssymb}
\usepackage{amsmath}
\usepackage{esint}
\usepackage{dcolumn}% Align table columns on decimal point
\usepackage{array}
\usepackage[letterpaper,%
left=0.75in,right=0.75in,top=1in,bottom=0.75in]{geometry}
\theoremstyle{plain}
\newtheorem{thm}{\protect\theoremname}
\theoremstyle{remark}
\newtheorem{rem}[thm]{\protect\remarkname}
\providecommand{\remarkname}{Remark}
\providecommand{\theoremname}{Theorem}
\usepackage{pgfplots}
\pgfplotsset{compat=newest,compat/show suggested version=false}
\usepackage{setspace}
\usepackage{tabu}	
\usepackage{upgreek}	
\usepackage{dsfont}
\usepackage{tensor}
\usepackage{fixmath}
\usepackage{nccmath}
\usepackage[bookmarks=false]{hyperref}
\hypersetup{
	colorlinks   = true,
	linkcolor = black,
	citecolor    = blue, %black,
	urlcolor = black,
	pdfstartview=%
}
\usepackage[numbers]{natbib}
\usepackage{etoolbox} % for '\AtBeginEnvironment' macro
\AtBeginEnvironment{pmatrix}{\everymath{\displaystyle}}
\AtBeginEnvironment{bmatrix}{\everymath{\displaystyle}}
\AtBeginEnvironment{array}{\everymath{\displaystyle}}
\DeclareMathOperator{\trace}{tr}

\newcommand{\aweak}[2]{\left\langle #1,\,#2 \right\rangle}
\newcommand{\tweak}[2]{\langle #1,\,#2 \rangle}

\newcommand{\vardif}[2]{\frac{\delta\mspace{-2.4mu} #1}{\delta\mspace{-1.6mu} #2}}
\newcommand{\widevardif}[2]{\frac{\delta\mspace{-1.2mu} #1}{\delta\mspace{-1.6mu} #2}}

\newcommand{\minivardif}[2]{\tfrac{\delta\mspace{-2.4mu} #1}{\delta\mspace{-1.6mu} #2}}
\newcommand{\wideminivardif}[2]{\tfrac{\delta\mspace{-1.2mu} #1}{\delta\mspace{-1.6mu} #2}}
\newcommand{\total}[2]{\frac{\mathrm{d}#1}{\mathrm{d#2}}}

\newcommand{\grad}{\nabla}
\newcommand{\diver}{\mathop{\mathrm{div}}\nolimits} % divergence

\newcommand{\tpd}[2]{\frac{\partial #1}{\partial #2}}
\newcommand{\dx}{\mathrm{d}x}
\newcommand{\dS}{\mathrm{d}A}
\newcommand{\op}[3]{\mathcal{#1}^{\hspace*{.2mm}{\raisebox{.6mm}{\scalebox{0.6}[0.6]{#2}}}}_{{#3}}}

\newcommand{\rboundbrak}[3]{	\left[#1,#2\right]_{\text{bound}}#3}

\newcommand{\boundbrak}[3]{	\left\{#1,#2\right\}_{\text{bound}}#3}

\newcommand{\dcont}{\raisebox{.3mm}{$:$}}
\newcommand{\contr}{\raisebox{.3mm}{$\cdot$}}
\def\zhet{{\boldsymbol{z}}}

\def\Om{{\mathrm{\Omega}}}
\def\opdot{\left[\,\cdot\,\right]}

\def\E{\mathbb{E}}
\def\R{\mathbb{R}}
\author{A. Moses Badlyan $^{1}$, B. Maschke  $^{2}$, C. Beattie  $^{3}$, V. Mehrmann $^{1}$% <-this % stops a space
	\thanks{$^{1}$ Institut f\"ur Mathematik, Sekr. MA 4-5, TU Berlin, Strasse des 17. Juni
		136, 10623 Berlin, Germany. \textit{E-mail}:~{badlyan@math.tu-berlin.de} (A. Moses Badlyan), {mehrmann@math.tu-berlin.de} (V. Mehrmann)}
	\thanks{$^{2}$	Univ Lyon, Universit{\'e} Claude Bernard Lyon 1, CNRS, LAGEP UMR 5007,
			Villeurbanne, 43 boulevard du 11 novembre 1918, F-69100, Villeurbanne, France. \textit{E-mail}:~{bernhard.maschke@univ-lyon1.fr}}
	\thanks{$^{3}$	Department of Mathematics,Virginia Tech, Blacksburg, VA 24061 USA. \textit{E-mail}:~{ beattie@vt.edu}}
	\thanks{The second author acknowledges the support of the Agence Nationale de la Recherche -- Deutsche Forschungsgemeinschaft (ANR-DFG), project INFIDHEM, ID ANR-16-CE92-0028. The first author would like to thank Christoph Zimmer (TU Berlin) for helpful conversations and comments.}
}

\begin{document}

\title{Open physical systems: from GENERIC 
	to port-Hamiltonian systems}
\maketitle

\begin{abstract}
Formulations of open physical systems within the framework of Non-Equilibrium Reversible/Irreversible Coupling (associated with the acronym ``GENERIC'') is related in this work with state-space realizations that are given as boundary port-Hamiltonian systems. This reformulation is carried out explicitly by splitting the dynamics of the system into a reversible contribution given by a Poisson bracket and an irreversible contribution given by a symmetric dissipation bracket, and is facilitated by the introduction of an exergy-like potential.\\
\end{abstract}


\section{Introduction}
For \emph{open} physical systems not subject to  dissipation, the extension of their Hamiltonian formulation in order to encompass port variables 
defined at the system boundary has led to
the introduction of \emph{port-Hamiltonian systems} defined with respect to
Stokes-Dirac structures~\cite{LeGorrecSIAM05,schaftGeomPhys02}.
For systems that include dissipative effects, boundary port-Hamiltonian realizations
have been defined using an implicit formulation of dissipation
\cite{Baaiu09b,CRC15_GreenProc,SCL_16}.

Another approach combining Hamiltonian and gradient system descriptions involves introducing a double-bracket two-generator formalism and is based on the same geometric structure as \emph{metriplectic systems}~\cite{BMR13}. This formalism has been introduced in~\cite{GrmOet97I, OetGrm97II}, under the acronym GENERIC (General Equations for the Non-Equilibrium Reversible/Irreversible Coupling).
In this formalism, the dynamics is split into a reversible contribution given by a Poisson bracket and an irreversible contribution given by a symmetric dissipation bracket. In~\cite{Oettinger_PhysRevE_06} the GENERIC formalism was extended to \emph{open} nonequilibrium thermodynamic systems, complementing the Poisson and dissipative brackets with boundary brackets.

In this paper, we relate these boundary brackets with pairs of conjugated port variables and show how the GENERIC framework for open physical systems may be reformulated into a port-Hamiltonian framework.
 
The paper is organized as follows. In Section~\ref{sec:fieldeq} we present the field equations.  
A state-space formulation of the dynamical system is given in Section~\ref{sec:SSF}.  A port-Hamiltonian reformulation is then presented in Section~\ref{sec:phformulation}.

\section{Field- and Constitutive Equations}\label{sec:fieldeq}
The considerations in this work are restricted to standard~\emph{Newtonian spacetime}~\cite{Rod95} of classical mechanics.
We denote by $\E^d$ the $d$-dimensional Euclidean manifold~\cite{ZeiQFT3}, a real $d$-dimensional Riemannian manifold whose elements are points and whose tangent spaces are isomorphic to the $d$-dimensional Hilbert space~${E}$, called~\emph{Euclidean space}. We assume that the~\emph{spatial manifold} at each instant of time is endowed with the geometric structure of the three-dimensional Euclidean manifold. Whenever convenient, the Euclidean manifold $\E^d$ is identified with $\R^d$.

The spatial domain $\Om$ is a non-empty open bounded and connected set in $\R^d$ with a boundary $\partial\mathrm{\Om}$ that is regular enough to guarantee the validity of the divergence theorem. The extensive quantities mass, linear momentum, and (internal) energy are monitored through the fields of their densities. Their time evolution and spatial distribution is described by the system of differential balance equations
\begin{subequations}\label{eq:fieldeqn}
	\begin{align}	
	%%%%%% evolution equation for mass density field
	\label{bal:mass}
	&\partial_t\rho
	+\diver\left(\rho{v}\right)=0,\\
	%%%%%% evolution equation for linear momentum density field
	\label{bal:mom}
	&\partial_t(\rho{v})
	+\diver\left(\rho{v}\otimes{v}\right) =
	\diver\mathrm{T},\\
	%%%%%% evolution equation for internal energy density field
	\label{bal:intenergy}
	&\partial_t(\rho \epsilon) + \diver(\rho \epsilon{v}) =
	-\diver{q}
	+\mathrm{T}\,\dcont\,\nabla{v},
	\end{align}
\end{subequations}
where $\rho$ is the field of the mass-density, ${M}=\rho{v}$ is the field of the linear momentum density, and $u=\rho\epsilon$ is the field of the internal energy density.
The system of partial differential equations~\eqref{eq:fieldeqn} is not a closed system and therefore must be supplemented by material specific~\emph{closure relations} called~\emph{constitutive equations}. In this work we restrict our considerations to linear irreversible thermodynamics and assume that the local equilibrium assumption of~\emph{classical irreversible thermodynamics} holds~\cite{Leb08}. Since we are interested in a Navier-Stokes-Fourier fluid we choose the following closure relations for the stress tensor $\mathrm{T}$, and the heat flux vector $q$, formulated for the three-dimensional $(d=3)$ domain~\cite{Mue85}
\begin{subequations}\label{eq:constrel}
	\begin{alignat}{4}
	\mathrm{T}&=-\,\mathrm{p}\,\mathrm{I}
	+(\zeta-\tfrac{2}{3}\eta)\diver({v})\mathrm{I}
	+\eta(\nabla{v}+\nabla{v}^\top),\label{cr:CauchyStr} \\
	{q}&=-\kappa\nabla\theta.\label{eq:Fourier}
	\end{alignat}
\end{subequations}
The second order tensor $\mathrm{T}$ \eqref{cr:CauchyStr} is called~\emph{Cauchy stress tensor} and takes values in $\mathbb{R}^{d\times d}_{\text{sym}}$, and the \emph{spatial velocity gradient}~$(\nabla{v})$ takes values in $\mathbb{R}^{d\times d}$. The coefficients $\zeta$ and $\eta$ are bulk and dynamic viscosity, respectively, $\mathrm{p}$ is the thermodynamic equilibrium pressure, $\theta$ is the absolute temperature field, and the coefficient $\kappa$ is the heat conductivity. The coefficients $\eta$, $\zeta$, and $\kappa$ will in general depend upon the temperature $\theta$ and the mass density $\rho$ and are defined to be nonnegative, $\kappa,\eta,\zeta\geq 0$~\cite{Mue85}.
We denote by $\upsigma$ the~\emph{viscosity part} of the stress tensor~\eqref{cr:CauchyStr} given by
\begin{align}\label{cr:viscostress}
\upsigma:=\lambda\diver\!\left(v\right)\mathrm{I}
+\eta(\grad{v}+(\nabla v)^\top),
\end{align}
where we have defined $\lambda:=\zeta-\tfrac{2}{3}\eta$. %is the~\emph{first Lam\'e parameter}~\cite{Sal12}.
Then $\mathrm{T}=-\mathrm{p}\,\mathrm{I}+\upsigma$, where $\mathrm{I}$ denotes the identity \cite[p.~342]{AbrMarRat88}. 
  \par
We view the field equations~\eqref{eq:fieldeqn} as describing state-space dynamics, taking then as macroscopic state variable the associated tuple of fields, viz. $(\rho,{M}, u)$.
\begin{rem}
The operation denoted by a colon in the term  $\mathrm{T}\,\dcont\,\grad{v}$ appearing on the right-hand side of~\eqref{bal:intenergy} is called~\emph{double contraction} such that $\mathrm{T}\,\dcont\,\grad{v}\equiv\trace(\mathrm{T}\cdot(\grad v)^\top)$, where $\trace$ denotes the \emph{trace}, defined as the contraction of a second order tensor, cf.~\cite[p.~69]{MarHug94}. 
\end{rem}

\section{State Space Formulation}\label{sec:SSF}
In order to connect the extended GENERIC framework for open systems with the notion of dynamical systems viewed in a system theoretic sense, we reinterpret the time-evolution equation formulated in~\cite[Eq.~19]{Oettinger_PhysRevE_06}.  Let $H$ and $S$ be two real-valued functionals defined on the state space and representing the total energy and entropy, respectively, of the system. Suppose $A$ is an arbitrary but fixed real-valued state-dependent functional defined on the state space.  The resulting rewritten \emph{bracket formalism} is given by the following evolution equation
\begin{align}\label{eq:timeEvoBracket}
\total{A}{t}&=\left\{A,H\right\}+\left[A,S\right]-\boundbrak{A}{H}{}-\rboundbrak{A}{S}{},
\end{align}
where $\{\cdot,\cdot\}$ is the full Poisson bracket, and $[\cdot,\cdot]$ is the full dissipation bracket of the GENERIC formalism~\cite{Oettinger_PhysRevE_06}. For dynamical systems confined to a time-independent spatial domain $\Om$ with boundary $\partial\Om$ these two brackets are split into bulk and boundary contributions, viz.
\begin{alignat}{4}\label{eq:fullbrack}
\{A,B\}=\{A,B\}_{\text{bulk}}+\{A,B\}_{\text{bound}},
\end{alignat}
with an analogous splitting for the dissipation bracket $[\cdot,\cdot]$. Note that the dynamics is described by the bulk related brackets, i.e. in case of the Hamiltonian part
\begin{alignat}{4}
\{A,H\}_{\text{bulk}}&=\left\{A,H\right\}-\boundbrak{A}{H}{},
\end{alignat}
and there is a similar equation for the dissipation bracket. However, only the full brackets~\eqref{eq:fullbrack} have all properties known from the GENERIC formalism for isolated systems, i.e., constituting a Poisson and a dissipation bracket, respectively.

The right-hand side of~\eqref{eq:timeEvoBracket} written by means of functional derivatives becomes
\begin{align}\label{eq:timeEvo}
\frac{\mathrm{d}A}{\mathrm{d}t}=&\int_\Om\vardif{\mspace{-0.5mu}A}{\zhet}\!\cdot\!\left(\mathfrak{J}\vardif{\mspace{-0.5mu}H}{\zhet}+\mathfrak{R}\vardif{S}{\zhet}\right)\dx\notag\\[5pt]
&-\int_{\partial\Om}\vardif{\mspace{-0.5mu}A}{\zhet}\!\cdot\!\left(\op{\mathfrak{J}}{$\partial$}{}\vardif{\mspace{-0.5mu}H}{\zhet}+\op{\mathfrak{R}}{\!$\partial$}{}\vardif{S}{\zhet}\right)\dS.
\end{align}
The (local) functional derivative $\minivardif{A}{\mspace{0.5mu}\zhet}$ of a real-valued functional $A$ is also known as the~\emph{(local) Volterra variational derivative}. For details and its relation to the~{G\^{a}teaux} and {Fr\'{e}chet} differentials see e.g.~\cite[p.~103~f.]{AbrMarRat88}, and~\cite{Ham82}. %~\cite{AveSmo67,Ham82}.

The reinterpretation of the GENERIC framework for open systems leads to the time evolution equation~\eqref{eq:timeEvo}. Based on this observation we define a system of operator equations, which combines a Hamiltonian and a gradient system that interact with the environment in a system theoretic sense through boundary ports.
All problems considered in this work are special cases of an abstract dynamical system represented by this system of operator equations, which has the following general form (for details see~\cite{MosZ17}):
\begin{subequations}\label{eq:Opeq}
	\begin{alignat}{3}
	\dot{\boldsymbol{z}} &= \,\mathcal{J}\mspace{-5mu}\left(\zhet\right)\!\vardif{H}{\boldsymbol{z}}\,\, +\,\, &&\mathcal{R}\!\left(\zhet\right)\!\vardif{S}{\boldsymbol{z}}\,+\, \mathcal{B}(\zhet)\mathbf{u} \quad &&\text{ in } \mathcal{D}_\zhet^\ast,
	\label{eq:operator_equation_open_I}\\[5pt]
	\!\!\mathbf{y}_{H} &=\mathcal{B}^\ast\!\!\left(\zhet\right)\!\vardif{H}{\zhet} && &&\text{ in } \mathcal{D}_{\mathbf{u}}^\ast,\label{eq:operator_equation_open_II}\\[5pt]
	\mathbf{y}_{S} &= &&\mathcal{B}^\ast\!\!\left(\zhet\right)\!\vardif{S}{\zhet} &&\text{ in } \mathcal{D}_{\mathbf{u}}^\ast. \label{eq:operator_equation_open_III}
	\end{alignat}
\end{subequations}
The space $\mathcal{D}_\zhet$ is a reflexive Banach space consisting of functions which map the domain $\Om$ into $\mathbb{R}^N$.
The state $\zhet$ evolves in the state space $\mathcal{Z}$, which in this work is an open subset of $\mathcal{D}_\zhet$, i.e. $\zhet\colon\mathbb{I}\rightarrow\mathcal{Z}$, where the time interval $\mathbb{I}$ is a bounded interval of the real line. Let $H,S\in C^\infty(\mathcal{Z})$ be the system\textquotesingle s total energy and entropy, respectively.
The linear \emph{bounded} operators $\mathcal{J}\!(\zhet)[\,\cdot\,],\mathcal{R}(\zhet)[\,\cdot\,]\colon \mathcal{D}_\zhet\to \mathcal{D}^\ast_\zhet$ are related to the full Poisson operator $\mathfrak{J}(\zhet)$ and the full dissipation operator $\mathfrak{R}(\zhet)$ of the extended GENERIC formalism for open systems, which are linear spatial-differential operators~\cite{MosZ17}, see also~\cite{GrmOet97I,OetGrm97II}.
We call the function $\mathbf{u}$ \emph{the combined input} and assume that it maps the time interval into a reflexive space, $\mathbf{u}\colon\mathbb{I}\to\mathcal{D}_\mathbf{u}$.
The operator $\mathcal{B}^{}(\zhet)[\,\cdot\,]\colon \mathcal{D}_{\mathbf{u}} \to \mathcal{D}_\zhet^\ast$  is related to both boundary operators of the GENERIC formalism for open systems, and $\mathcal{B}(\zhet)\mathbf{u} \in \mathcal{D}_\zhet^\ast$ combines the boundary contributions of the Hamiltonian and the gradient system. The combined input port $\mathbf{u}$ is complemented by two~\emph{output} ports $\mathbf{y}_H$ and~$\mathbf{y}_S$, defined by the relations \eqref{eq:operator_equation_open_II} and \eqref{eq:operator_equation_open_III}.
Assume that the total functional derivatives $\minivardif{H}{\zhet}$ and $\minivardif{S}{\zhet}$ exist and are elements of $\mathcal{D}_\zhet$. The operators $\mathcal{J}\!(\zhet)$ and $\mathcal{R}(\zhet)$ have to be modeled such that the following two \emph{non-interacting conditions} are satisfied
\begin{alignat}{6}\label{eq:NIC}
	&\mathcal{J}\!\left(\zhet\right)\!\vardif{S}{\zhet}=0 &\quad&\text{and}&\quad& \mathcal{R}\!\left(\zhet\right)\!\vardif{H}{\zhet}= 0\,.
\end{alignat}
Typically, the function spaces are chosen explicitly after the operators appearing in \eqref{eq:Opeq} are specified for the problem under consideration.  In the following, the case of classical hydrodynamics is considered;  the dynamical system represented by the field equations~\eqref{eq:fieldeqn} is supplemented by the linear closure relations \eqref{eq:constrel}, cf.~\cite[Sec.~III]{Oettinger_PhysRevE_06}.

\subsection{Classical Hydrodynamics}\label{sec:clHydrodyn}
Let the spatial domain $\Om\subset\mathbb{R}^3$ be a nonempty open bounded and connected set with Lipschitz boundary~\cite[p.~232]{Zei86} that does not change in time.
Let $W^{1,p}(\Om)$, $1\leq p<\infty$, denote the Sobolev space of all real valued functions with weak derivative where the function and its derivative are measurable and integrable up to the power of $p$ and hence elements of $L^p(\Om)$. We denote by $W^{1,p}(\Om)^\ast$ the dual of $W^{1,p}(\Om)$, for details see~\cite{AdaFou03}. We write $W^{1,p}(\Om;\mathbb{R}^N)$ for the space of mappings $f$ on $\Om$ with values in $\mathbb{R}^N$ such that each component function of $f$ is in $W^{1,p}(\Om)$. We write $W^{1,p}(\Om;\mathbb{R}^N)^\ast$ for its dual space.
We assume that $\mathcal{D}_\zhet$ is a closed subspace of $W^{1,p}(\Om;\mathbb{R}^N)$.
The state variable of classical hydrodynamics is an abstract function of the form
\begin{equation}%
\zhet = \begin{bmatrix}
\rho & {M} & u
\end{bmatrix}^\top\,,
\end{equation}%
that maps the time interval $\mathbb{I}$ into $\mathcal{Z}\subset\mathcal{D}_\zhet:= W^{1,3}(\Om;\R^N)$, 
where $\mathcal{Z}=\left\{\zhet\in\mathcal{D}_\zhet\,\vert\,\rho\geq \delta \text{ a.e. for some } \delta>0\right\}$ is an open subset of $\mathcal{D}_\zhet$, and $\zhet_t = [z_1, \ldots, z_N]^\top$, with $N=5$. Note that in classical hydrodynamics the internal energy density $u$ is amongst the independent state variables and the entropy density $s$ is the thermodynamic potential field~\cite[Sec.~II]{OetGrm97II}.

We have chosen~$W^{1,3}(\Om)$, since the component-functions of ${v}={M} / \rho $ and those of its derivative given by $$\nabla{v}=(\rho \nabla {M} -{M} \otimes \nabla \rho)/ \rho^2,$$ are element of $L^{p}(\Om)$ for all $p\in [1,3)$ by the continuous embedding of $W^{1,3}(\Om)$ into $L^q(\Om)$, $1\leq q<\infty$, \cite[Thm.~4.12]{AdaFou03}. Therefore the scalars ${v} \cdot {v}$, ${M} \cdot {v}$, ${M} \cdot {M}$ and the components of ${v}$ are also $W^{1,3}$-functions if one assumes the mass density $\rho$ and the linear momentum density ${M}$ to be slightly more regular, for details see \cite{MosZ17}.

The Hamiltonian $H$ is given by the physical energy function, i.e. the sum of kinetic and internal energy
\begin{align}\label{eq:total_energy}
	\begin{split}
	&H(\zhet)=\int_\Om h(\rho,M,u)\,\dx,\\
	&h(\rho,M,u):=\frac{{M}\cdot{M}}{2 \rho} + u,
	\end{split}
\end{align}
and the entropy $S$ (thermodynamic potential) is given by
\begin{equation}
\label{eq:entropy}
S(\zhet) = \int_\Om s(\rho,u)\, \dx.
\end{equation}
The functional derivatives of $H$ and $S$~\cite{Oettinger_PhysRevE_06} are
\begin{alignat}{3}\label{eq:vardifsHS}
&\vardif{H}{\zhet} = \biggr[\begin{array}{ccc}
\!\!\!-\frac{{v}\cdot{v}}{2} & {v}  & 1\!\!\!
\end{array}\biggr]^\top&\text{and}~~&	
\vardif{S}{\zhet} = \begin{bmatrix}
-\frac{\mu}{\theta}&{0}& \frac{1}{\theta}
\end{bmatrix}^\top\!\!\!.
\end{alignat}
We assume that the functional derivatives $\minivardif{H}{\zhet}$ and $\minivardif{S}{\zhet}$ are also elements of $\mathcal{D}_\zhet$. The chemical potential $\mu$ is related with the pressure via the thermodynamic constitutive relation
\begin{align}\label{eq:thermoconst}
\mathrm{p}+u = \theta s +\rho \mu,
\end{align}
that holds under the local equilibrium assumption~\cite{BerE94}.
We assume that for smooth enough state variable $\zhet\in \mathcal{Z}$ both operators $\mathcal{J}(\zhet)[\,\cdot\,], \mathcal{R}(\zhet)[\,\cdot\,] \colon W^{1,3}(\Om;\mathbb{R}^5)\to W^{1,3}(\Om;\mathbb{R}^5)^\ast$ are continuous.
The operator associated to the Hamiltonian part of the dynamics has the form
\begin{equation}\label{eq:JHydrodyn}
	\begin{split}
	\mathcal{J}\mspace{-6mu}\left(\zhet\right)=
	\left[\begin{array}{ccc}\displaystyle
	\begin{matrix}\\[-3ex]\!\!0\\[0.1ex]\end{matrix} &
	\begin{matrix}\\[-3ex]\!\!\op{\mathcal{J}}{}{\rho,{M}}\\[0.1ex]\end{matrix} &
	\begin{matrix}\\[-3ex]\!\!0\\[0.1ex]\end{matrix}
	\\[.5ex]%\hline
	\begin{matrix}\\[-3ex]\!\!\op{\mathcal{J}}{}{{M},\rho}\\[0.1ex]\end{matrix} &
	\begin{matrix}\\[-3ex]\!\!\op{\mathcal{J}}{}{{M},{M}}\\[0.1ex]\end{matrix} &
	\begin{matrix}\\[-3ex]\!\!\op{\mathcal{J}}{}{{M},u}\\[0.1ex]\end{matrix}
	\\[.5ex]%\hline
	\begin{matrix}\\[-3ex]\!\!0\\[0.1ex]\end{matrix} &
	\begin{matrix}\\[-3ex]\!\!\op{\mathcal{J}}{}{u,{M}}\\[0.1ex]\end{matrix} &
	\begin{matrix}\\[-3ex]\!\!0\\[0.1ex]\end{matrix}
	\end{array}\!\!\!\!\!\!\right].
	\end{split}
\end{equation}
In the weak formulation, the operators occupying the entries of $\mathcal{J}\!(\zhet)$ in~\eqref{eq:JHydrodyn} are defined implicitly through the relations
\begin{alignat}{2}
&\!\!\!\!\!\tweak{\varphi\indices{_\rho}}{\!\!\op{J}{}{\rho, {M}}\psi\indices{_{{M}}}}
=-\tweak{\psi\indices{_{{M}}}}{\!\!\op{J}{}{{M},\rho} \varphi\indices{_\rho} }
=\!\!\int_\Om\! \rho (\psi\indices{_{{M}}}\! \cdot\!\nabla) \varphi\indices{_\rho}\dx,
\end{alignat}
\begin{alignat}{2}
&\!\!\!\!\!\!\tweak{\varphi\indices{_{{M}}} }{\!\!\op{\mathcal{J}}{}{{M},{M}}\psi\indices{_{{M}}}}
=-\tweak{\psi\indices{_{{M}}}}{\!\!\op{\mathcal{J}}{}{{M},{M}}\varphi\indices{_{{M}}}}\notag\\
&\hspace*{1cm}= \int_\Om {M} \cdot \Big[(\psi\indices{_{{M}}}\cdot\nabla)\varphi\indices{_{{M}}} - (\varphi\indices{_{{M}}}\cdot\nabla)\psi\indices{_{{M}}}\Big]\dx,\\[10pt]
%%%%%%%%%
%
\label{eq:JuM}
&\!\!\!\!\!\tweak{\varphi\indices{_{u}} }{\!\!\op{J}{}{ u , {M}}\psi\indices{_{{M}}}} =  -  \tweak{\psi\indices{_{{M}}}}{\!\!\op{\mathcal{J}}{}{{M},u}\varphi\indices{_u} }\notag\\
&\hspace*{1cm}=\int_\Om u  (\psi\indices{_{{M}}} \cdot\nabla) \varphi\indices{_u} +
(\psi\indices{_{{M}}} \cdot\nabla) \left(\varphi\indices{_u}\,\mathrm{p}\right) \dx,
\end{alignat}
where ${\varphi}, {\psi} \in W^{1,3}(\Om;\mathbb{R}^5)$ are of the form $[\psi\indices{_\rho}~\psi\indices{_M}~\psi\indices{_u}]$, respectively. Note that the test function associated with the linear momentum density $M$ is vector valued, $\psi_M\colon\Om\to\mathbb{R}^3$.\par
The operator associated with the dissipative part of the dynamics has the form
\begin{equation}\label{eq:RHydrodyn}
	\begin{split}
	\mathcal{R}\mspace{-3mu}\left(\zhet\right)=
	\left[\begin{array}{ccc}\displaystyle
	\begin{matrix}\\[-3ex]0\\[0.1ex]\end{matrix} &
	\begin{matrix}\\[-3ex]0\\[0.1ex]\end{matrix} &
	\begin{matrix}\\[-3ex]0\\[0.1ex]\end{matrix}
	\\[.5ex]%\hline
	\begin{matrix}\\[-3ex]0\\[0.1ex]\end{matrix} &
	\begin{matrix}\\[-3ex]\op{\mathcal{R}}{}{{M,M}}\\[0.1ex]\end{matrix} &
	\begin{matrix}\\[-3ex]\op{\mathcal{R}}{}{{M,u}}\\[0.1ex]\end{matrix}
	\\[.5ex]%\hline
	\begin{matrix}\\[-3ex]0\\[0.1ex]\end{matrix} &
	\begin{matrix}\\[-3ex]\op{\mathcal{R}}{}{{u,M}}\\[0.1ex]\end{matrix} &
	\begin{matrix}\\[-3ex]\op{\mathcal{R}}{}{{u,u}}\\[0.1ex]\end{matrix}
	\end{array}\!\!\!\!\!\!\right].
	\end{split}
\end{equation}
The operators contained within $\mathcal{R}(\zhet)$ are defined through
\begin{alignat}{2}
\label{eq:RMM}
&\tweak{\varphi\indices{_{{M}}} }{\!\!\op{\mathcal{R}}{}{{M},{M}}\psi\indices{_{{M}}}} =\notag\\
&\begin{aligned}[b]\int_\Om\frac{\eta \theta}{2}
\Big[\grad\varphi\indices{_{{M}}} &+ \grad\varphi\indices{_{{M}}}{\!}^\top\Big]\dcont
	\Big[\grad\psi\indices{_{{M}}}+\grad\psi\indices{_{{M}}}{\!}^\top\Big]\\
&+\lambda\theta\diver\left(\varphi\indices{_{{M}}}\right)\diver\left(\psi\indices{_{{M}}}\right)\dx,\end{aligned}\\[10pt]
%%%%%%%%%%%%%%%%%%%%%
\label{eq:RMu}
&\tweak{\varphi\indices{_{{M}}} }{\!\!\op{\mathcal{R}}{}{{M},u}\psi\indices{_u}} =
\tweak{\psi\indices{_u} }{\!\!\op{\mathcal{R}}{}{u,{M}}\varphi\indices{_{{M}}}}=\notag\\
&\begin{aligned}[t]
\!\!\!\int_\Om \!
\!\!(-&\frac{\eta\theta}{2}\!\Big[\grad\varphi\indices{_{{M}}}
+\grad\varphi\indices{_{{M}}}{\!\!\!\!}^\top\Big]\dcont\,\mathrm{D}
-\!\frac{\lambda\theta}{2}\diver\!\left(\mspace{-1.2mu}\varphi\indices{_{{M}}}\mspace{-1.2mu}\right)\trace\!\left(\mathrm{D}\right))\psi\indices{_u}
\dx,\!\!\!\end{aligned}\!\!\!\!\!\!\\[10pt]
%%%%%%%%%%%%%%%%%%%%%%%%%%%%%%
\label{eq:Ruu}
%&\mbox{and} \notag\\
&\tweak{\varphi\indices{_u} }{\!\!\op{\mathcal{R}}{}{u,u}\psi\indices{_u}} =\notag\\
%&\begin{aligned}[t]
&\int_\Om\!\Big(\frac{\eta \theta}{2}
\mathrm{D}\,\dcont\,
\mathrm{D}
+\frac{\lambda\theta}{4}\trace\left(\mathrm{D}\right)^2\!\Big)\varphi\indices{_u}\psi\indices{_u}\! +\!\kappa \theta^2\grad\varphi\indices{_u}\!\cdot\!
\grad\psi\indices{_u}
\dx,\!\!
\end{alignat}
%%%%%%%%%%%%
where $\mathrm{D}=\grad {v}+(\nabla{v})^\top$ and $\trace\left(\mathrm{D}\right)=2\diver(v)$.

The operator $\mathcal{J}(\zhet)$ of~\eqref{eq:JHydrodyn} is an everywhere-defined bounded linear operator on a real, reflexive Banach space on which one may observe $\aweak{\xi}{\mathcal{J}(\zhet)\varphi}=-\aweak{\varphi}{\mathcal{J}(\zhet)\xi}$, i.e. the operator is skew-adjoint $\mathcal{J}^\ast=-\mathcal{J}$.
In a similar way, the dissipation operator $\mathcal{R}(\zhet)$ of~\eqref{eq:RHydrodyn} is self-adjoint,
i.e. $\mathcal{R}^\ast=\mathcal{R}$,  and moreover, $\aweak{\xi}{\mathcal{R}(\zhet)\xi}\geq 0$ holds for all
$\xi\in W^{1,3}(\Om;\mathbb{R}^5)$. That is, the dissipation operator $\mathcal{R}(\zhet)$ is both
self-adjoint  and semi-elliptic.

We now consider the boundary operators of~\eqref {eq:Opeq}, viz.
\begin{alignat}{4}
&\mathcal{B}^{}(\zhet)[\,\cdot\,]\colon \mathcal{D}_\mathbf{u} \to \mathcal{D}^\ast_\zhet, &~~\text{and}~~&\mathcal{B}^{\ast}\!(\zhet)[\,\cdot\,]\colon \mathcal{D}_\zhet\to \mathcal{D}^\ast_\mathbf{u},
\end{alignat}
and specify the associated function spaces. The reflexive space $\mathcal{D}_\mathbf{u}$ and its dual are chosen as $\mathcal{D}_\mathbf{u}:=L^q(\partial\Om;\mathbb{R}^5)$ and $\mathcal{D}^\ast_\mathbf{u}:= L^{\tilde{q}}(\partial\Om;\mathbb{R}^5)$, with $1< q,\tilde{q}<\infty$ and $1=\tfrac{1}{q}+\tfrac{1}{\tilde{q}}$.
For the concrete problem 
we define the linear operator $\mathcal{B}(\zhet)\opdot \colon L^{2}(\partial\Om;\mathbb{R}^5) \to W^{1,3}(\Om;\mathbb{R}^5)^\ast$
by means of $\langle\cdot,\cdot\rangle\colon W^{1,3}(\Om;\mathbb{R}^5)\times W^{1,3}(\Om;\mathbb{R}^5)^\ast\to \mathbb{R}$
as the pairing
\begin{align}\label{eq:bpair}
\tweak{\varphi}{\!\mathcal{B}(\zhet)\mathbf{u}}:=
\int_{\partial \Om}\!\!\!&-\Big[\varphi\indices{_\rho}\rho \mathrm{u}\indices{_1} -\varphi\indices{_{{M}}}\,\cdot\,\Big(\mathrm{u}\indices{_{[3:5]}} - {M}\mathrm{u}\indices{_1}\Big)\notag\\ &+\varphi\indices{_u}\left(\Big[u+\mathrm{p}\Big]\mathrm{u}\indices{_1}+\mathrm{u}\indices{_2}\right)\Big]\dS,
\end{align}
where the combined input port, $\mathbf{u}\colon\mathbb{I}\to L^{2}(\partial\Om;\mathbb{R}^5)$, for a fixed time parameter has the form $\mathbf{u}_t= [\mathrm{u}\indices{_1}, \mathrm{u}\indices{_2}, \mathrm{u}\indices{_{[3:5]}}]^\top$. It has components specified by the following block vector
\begin{align}\label{eq:input}
\mathbf{u}&=\Big[\begin{array}{lcr}
\!\!{v}|_{\partial \Om}\, \contr\, \nu ~~~ & q|_{\partial \Om}\, \contr\,  \nu &~~~ \upsigma|_{\partial \Om}\, \contr\, \nu\!\!
\end{array}\Big]^\top,
\end{align}
where $\nu$ is the function representing the unit normal vector.
For smooth enough $\zhet$ the operator $\mathcal{B}(\zhet)[\,\cdot\,]$ is assumed to be continuous. The adjoint operator $\mathcal{B}^{\ast}\!(\zhet)[\,\cdot\,]\colon W^{1,3}(\Om;\mathbb{R}^5) \to L^{2}(\partial\Om;\mathbb{R}^5)$ is defined as that linear operator satisfying
\begin{align}\label{eq:bpairing}
	\tweak{\boldsymbol{\varphi}}{{\mathcal{B}}(\zhet) \mathbf{u}} = \tweak{ \mathbf{u}}{{\mathcal{B}}^{\ast}\!(\zhet)\boldsymbol{\varphi}},
\end{align}
where the pairing on the right side of~\eqref{eq:bpairing} is the duality pairing $\langle\cdot,\cdot\rangle\colon L^q(\partial\Om;\mathbb{R}^5)\times L^{\tilde{q}}(\partial\Om;\mathbb{R}^5) \to \mathbb{R}$. Through these relations the output ports defined as $\mathbf{y}_H=\mathcal{B}^\ast\!(\zhet)\minivardif{H}{\zhet}$ and $\mathbf{y}_S=\mathcal{B}^\ast\!(\zhet)\minivardif{S}{\zhet}$ can be calculated. The output ports $\mathbf{y}_H$ and $\mathbf{y}_S$ for a fixed time take the form $\mathbf{y}_t = [\mathrm{y}\indices{_1}, \mathrm{y}\indices{_2}, \mathrm{y}\indices{_{[3:5]}}]^\top$. The output port related to the change of the Hamiltonian is
\begin{align}
%%%%%% output related to change of Hamiltonian
\!\!\mathbf{y}_H&=\left[\begin{array}{lcr}\!\!
-\left(\frac{ M\,\contr\, {M}}{2\rho} + \mu\rho + \theta s\right)|_{\partial \Om}~ & -1  &~ {v}|_{\partial \Om}
\end{array}\!\!\right]^\top\!\!\!\!,\label{eq:yh}
\intertext{where in \eqref{eq:yh} we have $\mu\rho + \theta s = u + \mathrm{p}$ \eqref{eq:thermoconst}. The output related to the change of the total entropy is given by}
%%%%% output related to change of total entropy
\mathbf{y}_S&=\left[\begin{array}{lcr}\!\!
- s|_{\partial \Om} ~~& -\frac{1}{\theta}|_{\partial \Om}  &~~~{0}\end{array}\right]^\top.\label{eq:ys}
\end{align}

\subsection{Balance of Energy and Entropy}
The mathematical model of classical hydrodynamics given as open infinite-dimensional nonlinear dissipative dynamical system in the operator setting~\eqref{eq:Opeq} satisfies the first and second law of thermodynamics, i.e., the balance of energy and the entropy inequality which for isolated systems reduce to the conservation of energy and production of entropy.
The time evolution of the Hamiltonian~\eqref{eq:total_energy} is given by
\begin{alignat}{4}\label{eq:totalH}
\!\!\!\!\!\frac{\mathrm{d}H}{\mathrm{d}t}
	&=\aweak{\vardif{H}{\zhet}}{\dot{\zhet}}
	 \overset{\eqref{eq:operator_equation_open_I}}{=}\aweak{\vardif{H}{\zhet}}{\mathcal{J}\vardif{H}{\zhet}+\mathcal{R}\vardif{S}{\zhet}+\mathcal{B}\mathbf{u}}\notag\\
	&=\aweak{\vardif{H}{\zhet}}{\mathcal{J}\vardif{H}{\zhet}}
		+\aweak{\vardif{H}{\zhet}}{\mathcal{R}\vardif{S}{\zhet}}
		+\aweak{\vardif{H}{\zhet}}{\mathcal{B}\mathbf{u}}\!,\!\!
\end{alignat}
where the time-derivative of the state variable is understood in the weak sense, i.e.  $\dot{\zhet}(t)\in\mathcal{D}_\zhet^\ast$ for almost every $t\in \mathbb{I}$ and $\|\dot{\zhet}\|_{\mathcal{D}_\zhet^\ast}$ is at least an element of $L_{\text{loc}}^{1}(\mathbb{I})$, see~\cite[Ch.~23.5]{Zei90}. Since the operator $\mathcal{J}(\zhet)$ is skew-adjoint, the balance equation~\eqref{eq:totalH} becomes
\begin{alignat}{4}\label{eq:totalH2}
&\frac{\mathrm{d}H}{\mathrm{d}t}
=\aweak{\vardif{S}{\zhet}}{\!\!\mathcal{R}\vardif{H}{\zhet}}
+\aweak{\vardif{H}{\zhet}}{\!\mathcal{B}\mathbf{u}}
\overset{\eqref{eq:NIC}}{=}\aweak{\vardif{H}{\zhet}}{\!\mathcal{B}\mathbf{u}}\notag\\[5pt]
=&-\int_{\partial\Om}\!\!\nu\cdot\left[\left(\frac{M\cdot M}{2\rho}+ u+\mathrm{p}\right){v}+q-(\upsigma^\top\!\!\cdot{v})\right]\!\!\dS,\!\!\!\!
\end{alignat}
where relations and definitions~\eqref{cr:viscostress},~\eqref{eq:vardifsHS},~\eqref{eq:thermoconst},~\eqref{eq:bpair},~\eqref{eq:input} and~\eqref{eq:yh} have been used. The balance equation~\eqref{eq:totalH2} corresponds to the integral total energy balance known from continuum physics and states that the energy of the system inside the spatial domain can only change due to convective and non-convective transport of energy over the boundary, cf. Equation (38) in~\cite{Oettinger_PhysRevE_06}. Note that influx terms like external body forces or thermal radiation are neglected.

Similarly, one shows that the total change of the entropy $S$ has the form of the entropy inequality (entropy balance) and corresponds to the second law of thermodynamics. Therefore, we formulate the time change of the functional $S$, given by~\eqref{eq:entropy}, and obtain
\begin{alignat}{4}\label{eq:totalS}
\frac{\mathrm{d}S}{\mathrm{d}t}
&=\aweak{\vardif{S}{\zhet}}{\dot{\zhet}}
\overset{\eqref{eq:operator_equation_open_I}}{=}\aweak{\vardif{S}{\zhet}}{\mathcal{J}\vardif{H}{\zhet}+\mathcal{R}\vardif{S}{\zhet}+\mathcal{B}\mathbf{u}}\notag\\
&=-\aweak{\vardif{H}{\zhet}}{\mathcal{J}\vardif{S}{\zhet}}
+\aweak{\vardif{S}{\zhet}}{\mathcal{R}\vardif{S}{\zhet}}
+\aweak{\vardif{S}{\zhet}}{\mathcal{B}\mathbf{u}}\notag\\
&\overset{\eqref{eq:NIC}}{=}\aweak{\vardif{S}{\zhet}}{\mathcal{R}\vardif{S}{\zhet}}
+\aweak{\vardif{S}{\zhet}}{\mathcal{B}\mathbf{u}}.
\end{alignat}
Assuming the operator $\mathcal{R}$ to be self-adjoint and semi-elliptic, we obtain a lower bound for the change of entropy
 \begin{alignat}{4}\label{ineq:totalS}
 \frac{\mathrm{d}S}{\mathrm{d}t}
 &=\aweak{\vardif{S}{\zhet}}{\mathcal{R}\vardif{S}{\zhet}}
 +\aweak{\vardif{S}{\zhet}}{\mathcal{B}\mathbf{u}}
 \geq\aweak{\vardif{S}{\zhet}}{\mathcal{B}\mathbf{u}}\notag\\
 &=-\int_{\partial \Om}\nu\cdot\left[s{v}+
 \frac{1}{\theta}q\right]\,\dS.
 \end{alignat}
Equation~\eqref{ineq:totalS} is the entropy balance for open systems, cf. Equation~(39) in~\cite{Oettinger_PhysRevE_06}.
The energy balance \eqref{eq:totalH2} and entropy balance \eqref{eq:totalS} expressed with the combined input $\mathbf{u}$ and the corresponding output ports $\mathbf{y}_H$ and $\mathbf{y}_S$ can be written as
\begin{alignat}{6}
	&\frac{\mathrm{d}H}{\mathrm{d}t}  = \tweak{\mathbf{y}_H}{\mathbf{u}} &\quad&\text{and}&\quad&
	\frac{\mathrm{d}S}{\mathrm{d}t} \geq \tweak{\mathbf{y}_S}{\mathbf{u}}.
\end{alignat}
%%%%%%%%%%%%%%%%%%%%%%%%%%%%

In this section the functionals and operators of \eqref{eq:Opeq} for classical hydrodynamics have been presented. The resulting state space framework encodes a weak formulation of the non-linear field equations \eqref{eq:fieldeqn} and complementary closure relations \eqref{eq:constrel}, as is shown by way of a simplified example in the upcoming section.

\section{Relation to port-Hamiltonian Systems}\label{sec:phformulation}
We proceed to use the mathematical model of classical hydrodynamics, as given in the operator setting~\eqref{eq:Opeq} that we specified in Section~\ref{sec:clHydrodyn}, in order to formulate a dissipative dynamical system that may be seen as representing a generalized port-Hamiltonian system. For the sake of simplicity, viscosity-induced dissipation is neglected.

\subsection{Exergy-like generating functional}
By means of the energy functional $H$ and entropy functional $S$ introduced in Section~\ref{sec:clHydrodyn} we define an \emph{exergy-like}~\cite{Alt17} energy functional $\mathcal{E}\in C^\infty(\mathcal{Z})$ by
\begin{alignat}{4}\label{eq:Exergy}
	&\mathcal{E}:= H - \mathcal{S} &~~\text{with}~~&\mathcal{S}:=\tau^{-1}_\circ S.
\end{alignat}
The functionals $H$ and $S$ are the physical energy function and the total entropy given by~\eqref{eq:total_energy} and \eqref{eq:entropy}, respectively. The term $\tau_\circ:=(1/\theta_\circ)\in\mathbb{R}_{\geq 0}$ is a scalar quantity representing a fixed reciprocal reference temperature value. Assume that for a smooth enough fixed state, $\zhet\in\mathcal{Z}$, the total functional derivative of the functional $\mathcal{E}$~\eqref{eq:Exergy} exists and is an element of $\mathcal{D}_\zhet$. Then it will be uniquely determined and satisfy
\begin{alignat}{6}\label{eq:ExergyTotvardif}
&\vardif{\mspace{1mu}\mathcal{E}}{\zhet}
=\vardif{H}{\zhet}-\vardif{\mspace{1mu}\mathcal{S}}{\zhet}.
\end{alignat}
By means of the partial functional derivatives of the Hamiltonian $H$ and the total entropy $S$, given in \eqref{eq:vardifsHS}, the partial functional derivatives of $\mathcal{E}$ can be calculated.

Now consider the dissipation operator $\mathcal{R}(\zhet)$~\eqref{eq:RHydrodyn} and its component operators defined via~\eqref{eq:RMM},~\eqref{eq:RMu}, and \eqref{eq:Ruu}. By setting the bulk and dynamic viscosity coefficients to zero, $\eta=\zeta=0$, the \emph{first Lam\'e parameter} vanishes identically, i.e. $\lambda=\zeta-\tfrac{2}{3}\eta=0$. This results in the simplified operator
\begin{equation}
	\begin{split}\label{op:Rsimple}
	\op{\mathcal{R}}{}{}\!(\zhet)=
	\left[\begin{array}{ccc}\displaystyle
	\begin{matrix}\\[-3ex]0\\[0.1ex]\end{matrix} &
	\begin{matrix}\\[-3ex]0\\[0.1ex]\end{matrix} &
	\begin{matrix}\\[-3ex]0\\[0.1ex]\end{matrix}
	\\[.5ex]%\hline
	\begin{matrix}\\[-3ex]0\\[0.1ex]\end{matrix} &
	\begin{matrix}\\[-3ex]0\\[0.1ex]\end{matrix} &
	\begin{matrix}\\[-3ex]0\\[0.1ex]\end{matrix}
	\\[.5ex]%\hline
	\begin{matrix}\\[-3ex]0\\[0.1ex]\end{matrix} &
	\begin{matrix}\\[-3ex]0\\[0.1ex]\end{matrix} &
	\begin{matrix}\\[-3ex]\op{\mathcal{R}}{}{u,u}\\[0.1ex]\end{matrix}
	\end{array}\!\!\!\!\!\!\right],
	\end{split}
\end{equation}
containing a single non-vanishing component-operator, $\mathcal{R}_{u,u}$, defined implicitly via
%%%%%%%%%%%% R_u,u
\begin{alignat}{2}\label{eq:RsimpleUU}
\tweak{\varphi\indices{_u} }{\!\!\op{\mathcal{R}}{}{u,u}\psi\indices{_u}} =
\int_\Om\kappa \theta^2\grad\varphi\indices{_u}\!\cdot\!
\grad\psi\indices{_u}
\dx.
\end{alignat}
%%%%%%%%%%%%
In~\cite{Ed98} it was observed that many dissipation operators $\mathfrak{R}$ of the original GENERIC framework can be expressed in the factorized form $\mathfrak{R}=\mathfrak{C}\mathfrak{D}\mathfrak{C}^\ast$. Motivated by this idea we multiply the dissipation operator~\eqref{op:Rsimple} by the scalar $\tau_\circ$ and assume that the resulting operator can be factorized as
\begin{align}\label{eq:factR}
	\tau_\circ\mathcal{R}=\mathcal{C}\mathcal{D}\mathcal{C}^\ast.
\end{align} The operator $\mathcal{C}(\zhet)$ at the right-hand side of~\eqref{eq:factR} has the form
\begin{equation}
\begin{split}\label{op:C}
\op{\mathcal{C}}{}{}\!(\zhet)=
\left[\begin{array}{ccc}\displaystyle
\begin{matrix}\\[-3ex]0\\[0.1ex]\end{matrix} &
\begin{matrix}\\[-3ex]0\\[0.1ex]\end{matrix} &
\begin{matrix}\\[-3ex]0\\[0.1ex]\end{matrix}
\\[.5ex]%\hline
\begin{matrix}\\[-3ex]0\\[0.1ex]\end{matrix} &
\begin{matrix}\\[-3ex]0\\[0.1ex]\end{matrix} &
\begin{matrix}\\[-3ex]0\\[0.1ex]\end{matrix}
\\[.5ex]%\hline
\begin{matrix}\\[-3ex]0\\[0.1ex]\end{matrix} &
\begin{matrix}\\[-3ex]0\\[0.1ex]\end{matrix} &
\begin{matrix}\\[-3ex]\op{\mathcal{C}}{}{u,u}\\[0.1ex]\end{matrix}
\end{array}\!\!\!\!\!\!\right],
\end{split}
\end{equation}
where $\op{\mathcal{C}}{}{u,u}\!\opdot\colon L^{\frac{3}{2}}(\Om,\mathbb{R}^d)\to W^{1,3}(\Om,\mathbb{R})^\ast$ is defined via
\begin{alignat}{2}
\tweak{\phi}{\op{\mathcal{C}}{}{u,u}\psi}
%\equiv\tweak{\phi}{\op{\mathcal{C}}{}{{\cdot}u}\psi}
=\int_\Om \grad\phi\cdot \psi\,\dx.
\end{alignat}
Significantly, the operator $\mathcal{C}(\zhet)$ is self-adjoint, $\mathcal{C}=\mathcal{C}^\ast$, i.e. one may observe $\tweak{\phi}{\op{\mathcal{C}}{}{u,u}\psi} = \tweak{\psi}{\op{\mathcal{C}}{}{u,u}\phi}$.
The operator $\mathcal{D}$ of the factorization~\eqref{eq:factR} is of the form
\begin{equation}
	\begin{split}\label{op:D}
	\op{\mathcal{D}}{}{}\!(\zhet)=
	\left[\begin{array}{ccc}\displaystyle
	\begin{matrix}\\[-3ex]0\\[0.1ex]\end{matrix} &
	\begin{matrix}\\[-3ex]0\\[0.1ex]\end{matrix} &
	\begin{matrix}\\[-3ex]0\\[0.1ex]\end{matrix}
	\\[.5ex]%\hline
	\begin{matrix}\\[-3ex]0\\[0.1ex]\end{matrix} &
	\begin{matrix}\\[-3ex]0\\[0.1ex]\end{matrix} &
	\begin{matrix}\\[-3ex]0\\[0.1ex]\end{matrix}
	\\[.5ex]%\hline
	\begin{matrix}\\[-3ex]0\\[0.1ex]\end{matrix} &
	\begin{matrix}\\[-3ex]0\\[0.1ex]\end{matrix} &
	\begin{matrix}\\[-3ex]\op{\mathcal{D}}{}{u,u}\\[0.1ex]\end{matrix}
	\end{array}\!\!\!\!\!\!\right],
	\end{split}
\end{equation}
where the operator $\op{\mathcal{D}}{}{u,u}\!\opdot\colon L^{\frac{3}{2}}(\Om,\mathbb{R}^d)^\ast\to  L^{\frac{3}{2}}(\Om,\mathbb{R}^d)$ is defined via
%%%%%%%%%%%%
\begin{alignat}{2}\label{eq:Duu}
\tweak{\xi}{\op{\mathcal{D}}{}{u,u}{\varphi}} =
\int_\Om \kappa \frac{\tau_\circ}{\tau^2}{\varphi}\cdot\xi\,\dx.
\end{alignat}
%%%%%%%%%%%%
One may observe that $\mathcal{D}(\zhet)$ is also self-adjoint, $\mathcal{D}=\mathcal{D}^\ast$.
Since the reciprocal absolute temperature $\tau=1/\theta$ takes nonnegative values and the heat-conductivity $\kappa$ is nonnegative, we have $(\tau_\circ/\tau^2)\kappa\geq 0$, such that $\aweak{\xi}{\mathcal{D}(\zhet)\xi}\geq 0$ for all $\xi\in L^{\frac{3}{2}}(\Om,\mathbb{R}^d)^\ast$. Therefore the operator $\mathcal{D}(\zhet)$ is both, self-adjoint and semi-elliptic.
The operator  $\mathcal{C}^\ast\!\!\left(\zhet\right)$ has the form
\begin{equation}
	\begin{split}\label{op:Cadj}
	\mathcal{C}^\ast\!(\zhet)=
	\left[\begin{array}{ccc}\displaystyle
	\begin{matrix}\\[-3ex]0\\[0.1ex]\end{matrix} &
	\begin{matrix}\\[-3ex]0\\[0.1ex]\end{matrix} &
	\begin{matrix}\\[-3ex]0\\[0.1ex]\end{matrix}
	\\[.5ex]%\hline
	\begin{matrix}\\[-3ex]0\\[0.1ex]\end{matrix} &
	\begin{matrix}\\[-3ex]0\\[0.1ex]\end{matrix} &
	\begin{matrix}\\[-3ex]0\\[0.1ex]\end{matrix}
	\\[.5ex]%\hline
	\begin{matrix}\\[-3ex]0\\[0.1ex]\end{matrix} &
	\begin{matrix}\\[-3ex]0\\[0.1ex]\end{matrix} &
	\begin{matrix}\\[-3ex]\mathcal{C}^\ast_{u,u}\\[0.1ex]\end{matrix}
	\end{array}\!\!\!\!\!\!\right],
	\end{split}
\end{equation}
where $\mathcal{C}^\ast_{u,u}\!\opdot\colon W^{1,3}(\Om,\mathbb{R})\to L^{\frac{3}{2}}(\Om,\mathbb{R}^d)^\ast$ is defined via
\begin{alignat}{2}\label{eq:Cadj}
\tweak{\tilde{\phi}}{\mathcal{C}^\ast_{u,u}\tilde{\psi}}
=\int_\Om \grad\tilde{\psi}\cdot \tilde{\phi}\,\dx.
\end{alignat}
Since according~\eqref{eq:vardifsHS} we have $\minivardif{H}{u}=1$, the operator $\mathcal{C}^\ast\!\!\left(\zhet\right)$ satisfies $\tweak{\tilde{\phi}}{\mathcal{C}^\ast\!\!\left(\zhet\right)\!\minivardif{H}{\zhet}} = 0$ for all $\tilde{\phi}\in L^{\frac{3}{2}}(\Om,\mathbb{R}^d)$. This leads to the following non-interacting condition, cf.~\cite[p.~68]{Oet05}
\begin{align}\label{NIC:Cadj}
\mathcal{C}^\ast\!\!\left(\zhet\right)\!\vardif{H}{\zhet}= 0.
\end{align}
Under the assumption that condition~\eqref{NIC:Cadj} holds, we obtain
\begin{alignat}{2}\label{NIC:CDC}
\mathcal{C}\mathcal{D}\mathcal{C}^\ast\!\!\left(\zhet\right)\!\widevardif{\mathcal{E}}{\zhet}&
\overset{\eqref{eq:ExergyTotvardif}}{=} -\mathcal{C}\mathcal{D}\mathcal{C}^\ast\!\!\left(\zhet\right)\!\widevardif{\mathcal{S}}{\zhet}\,.
\end{alignat}
 Furthermore, from the non-interacting condition~\eqref{eq:NIC} of the Poisson operator $\mathcal{J}$ we obtain
\begin{alignat}{2}\label{eq:NIC2}
\mathcal{J}(\zhet)\widevardif{\mathcal{E}}{\zhet}&=\mathcal{J}(\zhet)\vardif{H}{\zhet}.
\end{alignat}

Following \cite{HutS13} and restricted to linear irreversible thermodynamics, we assume that general thermodynamic forces, $F_{\text{orce}}$, and general thermodynamic fluxes, $F_{\text{lux}}$, are related through an operator equation having the form
\begin{align}\label{rel:ff}
	F_{\text{lux}}=\mathfrak{D}F_{\text{orce}}.
\end{align}
For heat-conduction, the left-hand side of Equation~\eqref{rel:ff} reflects the non-convective transport of internal energy, i.e. $F_{\text{lux}}\equiv q$. By choosing $F_{\text{orce}}=\nabla(\theta^{-1})$ and $\mathfrak{D}=\tau_\circ^{-1}\mathcal{D}_{u,u}$, where $\mathcal{D}_{u,u}$ is the operator~\eqref{op:D},  equation~\eqref{rel:ff} turns into the law of Fourier~\eqref{eq:Fourier},

We use relation~\eqref{rel:ff} as orientation and consequently define the dimensionless~\emph{thermodynamic force}
\begin{align}\label{eq:Fforce}
F_{\text{orce}}:=-\,\mathcal{C}^\ast\!\!\left(\zhet\right)\!\widevardif{\mathcal{E}}{\zhet}\overset{\eqref{NIC:Cadj}}{=}
\mathcal{C}^\ast\!\!\left(\zhet\right)\!\widevardif{\mathcal{S}}{\zhet}.
\end{align}
The right-hand side of~\eqref{eq:Fforce} may be thought of as the dimensionless analogue to the thermodynamic force, which drives the non-convective heat flux $q$. In accordance to equation~\eqref{rel:ff} we use~\eqref{eq:Fforce} to define the thermodynamic flux
\begin{align}\label{eq:Fflux}
F_{\text{lux}}:=-\mathcal{D}\mathcal{C}^\ast\!\!\left(\zhet\right)\!\widevardif{\mathcal{E}}{\zhet}\overset{\eqref{NIC:CDC}}{=}
\mathcal{D}\mathcal{C}^\ast\!\!\left(\zhet\right)\!\widevardif{\mathcal{S}}{\zhet}.
\end{align}
Note that although $F_\text{orce}$ as defined in \eqref{eq:Fforce} results in a dimensionless reciprocal~\emph{temperature} gradient given by $\left(\mathcal{C}^\ast_{u,u}\wideminivardif{\mathcal{S}}{u}\right)$, the corresponding flux $F_{\text{lux}}$ formulated according to~\eqref{eq:Fflux} is still the heat flux $q$. This is because of the scalar $\tau_\circ$ contained in the component-operator $\mathcal{D}_{u,u}$.

Next, we introduce one single output port, denoted by $\mathbf{y}_\mathcal{E}$, understood as the output of the dynamical system whose dynamics is generated by means of the exergy-like energy functional $\mathcal{E}$ \eqref{eq:Exergy}. For this we use the definitions of the output ports, $\mathbf{y}_H$ and $\mathbf{y}_S$, which were introduced in Section~\ref{sec:SSF} as part of the abstract operator representation of the extended GENERIC formalism for open systems~\eqref{eq:Opeq}. In accordance to \eqref{eq:operator_equation_open_II} and \eqref{eq:operator_equation_open_III} we define $\mathbf{y}_\mathcal{E}:=\mathcal{B}^\ast\!(\zhet)\wideminivardif{\mathcal{E}}{\zhet}$ such that
%%%%%%%%%%%%
\begin{align}\label{eq:ExergyInp}
\!\!\!\mathbf{y}_\mathcal{E}&=\mathcal{B}^\ast\!(\zhet)\widevardif{\mathcal{E}}{\zhet}
\overset{\eqref{eq:ExergyTotvardif}}{=}\mathcal{B}^\ast\!(\zhet)\!\left[\vardif{H}{\zhet}-\widevardif{\mathcal{S}}{\zhet}\right]\notag\\
&\overset{\eqref{eq:Exergy}}{=}\mathcal{B}^\ast\!(\zhet)\!\left[\vardif{H}{\zhet} -\tau^{-1}_\circ\vardif{S}{\zhet}\right]
=\mathbf{y}_H - \tau^{-1}_\circ \mathbf{y}_S\,.
\end{align}
%%%%%%%%%%%%
We define the output port $\mathbf{y}_\mathcal{S}:= \tau^{-1}_\circ \mathbf{y}_S$ and write~\eqref{eq:ExergyInp} as
\begin{align}\label{eq:ExergyInp2}
\mathbf{y}_\mathcal{E}
=\mathbf{y}_H - \mathbf{y}_\mathcal{S}.
\end{align}
\subsection{Port-Hamiltonian systems in input-output representations}
Let the energy functional $\mathcal{E}\in C^\infty(\mathcal{Z})$ be given by~\eqref{eq:Exergy} and the operators $\mathcal{J}(\zhet)[\,\cdot\,],(\mathcal{CDC}^\ast)(\zhet)[\,\cdot\,]\colon \mathcal{D}_\zhet\to \mathcal{D}^\ast_\zhet$ be defined as in \eqref{eq:JHydrodyn} and \eqref{op:Rsimple}, respectively. Assume that the \emph{combined input} function  $\mathbf{u}\colon\mathbb{I}\to\mathcal{D}_\mathbf{u}$, $t\mapsto \mathbf{u}_t= [\mathrm{u}\indices{_1}, \mathrm{u}\indices{_2}, \mathrm{u}\indices{_{[3:5]}}]^\top$ is identical with~\eqref{eq:input}, except that $\mathrm{u}\indices{_{[3:5]}}= 0$. Let the~\emph{output} function $\mathbf{y}_\mathcal{E}$ be defined through relation~\eqref{eq:ExergyInp}. Then the system of operator equations~\eqref{eq:Opeq} can be rewritten as
%%%%%%%%%%%%
\begin{subequations}\label{sys:pH1}
	\begin{alignat}{4}
	\label{eq:ExOp1}\!\!
	\dot{\zhet} &= \Bigr[\mathcal{J}(\zhet)-\mathcal{CDC}^\ast\!(\zhet)\Bigr]\widevardif{\mathcal{E}}{\zhet} + \mathcal{B}(\zhet)\mathbf{u}  &\quad\text{in}~\mathcal{D}_\zhet^\ast,&
	\\[0pt]
	\label{eq:ExOp2}\!\!
	\displaystyle\mathbf{y}_{\mathcal{E}} & =\mathcal{B}^\ast\!(\zhet) \widevardif{\mathcal{E}}{\zhet}
	&\quad\text{in}~\mathcal{D}_{\mathbf{u}}^\ast.&
	\end{alignat}
\end{subequations}
%%%%%%%%%%%%
The system~\eqref{sys:pH1} represents an infinite-dimensional nonlinear dissipative dynamical system given in the form of a combined Hamiltonian and a gradient system.
Note that the operator $\mathcal{J}(\zhet)$ in~\eqref{eq:ExOp1}, which is specified in~\eqref{eq:JHydrodyn}, via duality pairing induces a bracket that is identical with the full Poisson bracket of the GENERIC formulation of classical hydrodynamics, cf. Equation (29) in \cite{Oettinger_PhysRevE_06}. In other words, the bracket defined by $\left\{A,B\right\}\!\left(\zhet\right):=\aweak{\minivardif{A}{\zhet}}{\!\mathcal{J}(\zhet)\minivardif{B}{\zhet}}$ for arbitrary $A,B\in C^\infty(\mathcal{Z})$ and $\zhet\in\mathcal{Z}$, constitutes a Poisson bracket. Hence, it is anti-symmetric and satisfies both the Leibniz rule and the Jacobi identity.

The balance law for the exergy-like energy functional $\mathcal{E}$ \eqref{eq:Exergy} takes the form of a balance inequality. Under the assumption that the system under consideration is given as abstract dynamical system represented by the dissipative port-Hamiltonian system~\eqref{sys:pH1}, the following inequality holds
\begin{align}\label{eq:Exbal}
\frac{\mathrm{d}\mathcal{E}}{\mathrm{d}t}&=\aweak{\widevardif{\mathcal{E}}{\zhet}}{\dot{\zhet}}
\overset{\eqref{eq:ExOp1}}{=}\aweak{\widevardif{\mathcal{E}}{\zhet}}{\!\!\Big[\mathcal{J}-\mathcal{CDC}^\ast\Big]\widevardif{\mathcal{E}}{\zhet} + \mathcal{B}\mathbf{u}}\notag\\
&\hspace*{-.4em}\overset{\eqref{eq:ExOp2}}{=}\!\!-\aweak{\widevardif{\mathcal{E}}{\zhet}}{\mathcal{C}\mathcal{D}\mathcal{C}^\ast\widevardif{\mathcal{E}}{\zhet}}+\aweak{\mathbf{y}_\mathcal{E}}{\mathbf{u}}\leq\aweak{\mathbf{y}_\mathcal{E}}{\mathbf{u}} .
\end{align}
Similar to suggestions made in \cite{schaftGeomPhys02, SCL_16}, and used in the applications of
\cite{Baaiu09b,CRC15_GreenProc}, the system \eqref{sys:pH1} may be extended so as to give it the appearance of a skew-symmetric operator
\begin{align}\label{eq:dynmod1}
\left[\!\!\!\begin{array}{c}\displaystyle
\begin{matrix}\\[-2ex]\dot{\zhet}\\[1ex]\end{matrix}
\\[.5ex]%\hline
\begin{matrix}\\[-2ex]\,F_{\text{orce}}\\[1ex]\end{matrix}
\\[.5ex]%\hline
\begin{matrix}\\[-2ex]-\mathbf{y}_{\mathcal{E}}\\[1ex]\end{matrix}
\end{array}\!\!\!\!\right]
%---------------------
=\left[\begin{array}{ccc}
\begin{matrix}\\[-2ex]\mathcal{J}\mspace{-5mu}\left(\zhet\right)\\[1ex]\end{matrix} &
\begin{matrix}\\[-2ex]\mathcal{C}\!\left(\zhet\right)\\[1ex]\end{matrix} &
\begin{matrix}\\[-2ex]\mathcal{B}\!\left(\zhet\right)\\[1ex]\end{matrix}
\\[.5ex]%\hline
\begin{matrix}\\[-2ex]\!\!\!-\mathcal{C}^\ast\!\!\left(\zhet\right)\\[1ex]\end{matrix} &
\begin{matrix}\\[-2ex]0\\[1ex]\end{matrix} &
\begin{matrix}\\[-2ex]0\\[1ex]\end{matrix}
\\[.5ex]%\hline
\begin{matrix}\\[-2ex]\!\!\!-\mathcal{B}^\ast\!\!\left(\zhet\right)\\[1ex]\end{matrix} &
\begin{matrix}\\[-2ex]0\\[1ex]\end{matrix} &
\begin{matrix}\\[-2ex]0\\[1ex]\end{matrix}
\end{array}\!\!\!\!\right]
%---------------------	
\left[\!\!\!\begin{array}{c}\displaystyle
\begin{matrix}\\[-2ex]\delta_\zhet\mathcal{E}\\[1ex]\end{matrix}
\\[.5ex]%\hline
\begin{matrix}\\[-2ex]F_{\text{lux}}\\[1ex]\end{matrix}
\\[.5ex]%\hline
\begin{matrix}\\[-2ex]\mathbf{u}\\[1ex]\end{matrix}
\end{array}\!\!\!\!\!\right],
\end{align}
where $\delta_\zhet\mathcal{E}:=\widevardif{\mathcal{E}}{\zhet}$ and the relation \eqref{rel:ff} between the flux and the force variables has been used.
We suggest that the extended skew-adjoint operator in \eqref{eq:dynmod1} generates a Dirac structure on the product spaces of the Banach spaces on which the variables are defined, in the sense of \cite{IftSanGol05,Kurula_JMAA_10}. Thereby the system \eqref{eq:dynmod1} with   \eqref{rel:ff}  may be interpreted as a dissipative port-Hamiltonian system in descriptor form.
%%%%%%%%%%%%

\subsection{Heat conducting inviscid compressible fluid}\label{sec:Hydrodynamics}
Let $H$ be the Hamiltonian~\eqref{eq:total_energy} and $S$ the total entropy~\eqref{eq:entropy}.
Define the functional $\mathcal{S}(\zhet)=\int_\Om\tilde{s}(\rho,u)\dx$ with the density
\begin{align}
	\tilde{s}(\rho,u):=\tau_\circ^{-1} s(\rho,u),
\end{align} where $s(\rho,u)$ is the entropy density in~\eqref{eq:entropy}. Relation~\eqref{eq:thermoconst} expressed with density $\tilde{s}$ becomes
\begin{align}
	\mathrm{p}+u =\frac{\tau_\circ}{\tau}\tilde{s} +\rho \mu.
\end{align}
Using $H$ and $\mathcal{S}$, the energy functional $\mathcal{E}\in C^\infty(\mathcal{Z})$ can be formulated in accordance to~\eqref{eq:Exergy} such that
 \begin{align}
 &\mathcal{E}(\zhet)=\int_\Om e(\rho,M,u)\,\dx,\label{eq:functional}\\
 &e(\rho,M,u):=h(\rho,M,u)-\tilde{s}(\rho,u),
 \end{align}
 where the density $h(\rho,M,u)$ is specified in~\eqref{eq:total_energy}.
Then the total functional derivative of \eqref{eq:functional} is given by
\begin{alignat}{6}
\widevardif{\mathcal{E}}{\zhet}
&=\biggr[\!\!\begin{array}{lcr}
\widevardif{\mathcal{E}}{\rho} & \widevardif{\mathcal{E}}{M}   & \widevardif{\mathcal{E}}{u}   \end{array}\!\!\biggr]^\top\notag\\
&=\biggr[\!\!\begin{array}{lcr}
\Big(-\frac{v\cdot v}{2}+\frac{\tau}{\tau_\circ}\mu\Big)~~ & {v} &~~  \Big(1-\frac{\tau}{\tau_\circ}\Big)
\end{array}\!\!\biggr]^\top.\label{eq:Epvardifs}
\end{alignat}
The output port $\mathbf{y}_\mathcal{E}$ can be constructed according~\eqref{eq:ExergyInp2}.
The input port $\mathbf{u}$ in~\eqref{eq:ExOp1} is identical with~\eqref{eq:input}, a mapping of the form $\mathbf{u}\colon\mathbb{I}\to L^{2}(\partial\Om;\mathbb{R}^5)$, $t\mapsto \mathbf{u}_t= [\mathrm{u}\indices{_1}, \mathrm{u}\indices{_2}, \mathrm{u}\indices{_{[3:5]}}]^\top$. The components of the input and output ports are specified by
%%%%%%%%%%%%%%%%%%%%%%%%%%%
\begin{align}
\label{eq:inputPH}\!\!
\mathbf{u}&=\biggr[\begin{array}{lcr}
\!\!{v}|_{\partial \Om}\, \contr\, \nu ~~~ &  F_{\text{lux}}|_{\partial \Om}\, \contr\,  \nu &~~~ 0 \!\!
\end{array}\biggr]^\top,\\
\label{eq:outputPH}\!\!
\mathbf{y}_\mathcal{E}&\!=\!\begin{bmatrix}
-\!\left(\!\frac{M\,\contr\, {M}}{2\rho} + u + \mathrm{p} - \tilde{s}\!\right)\!|_{\partial \Om}\!\! &\! \frac{\tau }{\tau_\circ}|_{\partial \Om}-1\! &\! v|_{\partial \Om}
\end{bmatrix}^\top\!\!\!\!\!\!.\!\!\!
\end{align}~
%%%%%%%%%%%%%%%%%%%%%%%%%%%%
The boundary operator $\mathcal{B}^{}(\zhet)[\,\cdot\,]\colon \mathcal{D}_\mathbf{u} \to \mathcal{D}^\ast_\zhet$, and its adjoint~$\mathcal{B}^{\ast}\!(\zhet)[\,\cdot\,]\colon \mathcal{D}_\zhet\to \mathcal{D}^\ast_\mathbf{u}$ are identical with the corresponding boundary operators introduced in Section~\ref{sec:clHydrodyn}.
Since we have neglected viscosity-related dissipation, the pairing
specified in~\eqref{eq:bpair} can be written as
%%%%%%%%%%%%%%%%%%%%%%%%%%%%
\begin{align}\label{Bu:exergy}
&\tweak{\psi}{\!\mathcal{B}(\zhet)\mathbf{u}}=\notag\\
\!\!\!\!-&\begin{aligned}[4]
\int_{\partial \Om}\!\!\Big[\Big(\psi\indices{_\rho}\,\rho +\psi\indices{_{{M}}}\cdot{M}\Big)\mathrm{u}\indices{_1}\! +\psi\indices{_u}\!\left(\Big[u+\mathrm{p}\Big]\mathrm{u}\indices{_1}\!+\!\mathrm{u}\indices{_2}\right)\!\Big]\!\dS.
\end{aligned}\!\!\!
\end{align}
%%%%%%%%%%%%%%%%%%%%%%%%%%%%
The right-hand side of~\eqref{eq:Exbal} can be calculated using the input $\mathbf{u}$ and output $\mathbf{y}_\mathcal{E}$ as specified in~\eqref{eq:inputPH} and~\eqref{eq:outputPH}, respectively. From this we obtain
%%%%%%%%%%%%%%%%%%%%%%%%%%%%
\begin{align}\label{eq:exergybal}
&\aweak{\mathbf{y}_\mathcal{E}}{\mathbf{u}}=\notag\\
\!\!-&\!\!\begin{aligned}[t]
\int_{\partial\Om}\!\nu\cdot\left[\left(\frac{M\,\contr\, {M}}{2\rho} + u + \mathrm{p} - \tilde{s}\right)v -\frac{\tau-\tau_\circ}{\tau_\circ} F_\text{lux} \right]\!\dS.\!\!
\end{aligned}
\end{align}
%%%%%%%%%%%%%%%%%%%%%%%%%%%%
Equation~\eqref{eq:exergybal} is an upper bound for the change of the system total exergy-like energy $\mathcal{E}(\zhet)$. The first term in the integrand on the right-hand side of~\eqref{eq:exergybal} is the convective flux of ``exergy'', while the second term represents the non-convective flux, that  vanishes for $\tau=\tau_\circ$.
Note that the term
\begin{align}\label{eq:exergybal2}
	\aweak{\widevardif{\mathcal{E}}{\zhet}}{\mathcal{C}\mathcal{D}\mathcal{C}^\ast\widevardif{\mathcal{E}}{\zhet}}=\int_\Om  \kappa\frac{\tau_\circ}{\tau^2}\grad\!\left(\frac{\tau}{\tau_\circ}\right) \cdot\grad\!\left(\frac{\tau}{\tau_\circ}\right)\dx,
\end{align}
on the right-hand side of energy balance~\eqref{eq:Exbal} vanishes for $\tau_\circ=\tau$, and  with~\eqref{eq:exergybal} it follows that then balance equation~\eqref{eq:Exbal} takes the form 
%of the following balance equation
\begin{align}\label{eq:Exbal2}
\frac{\mathrm{d}\mathcal{E}}{\mathrm{d}t}=
-&\begin{aligned}[t]
\int_{\partial\Om}\nu\cdot\left[\left(\frac{M\,\contr\, {M}}{2\rho} + u + \mathrm{p} - \tilde{s}\right)v \right]\!\dS.
\end{aligned}
\end{align}
The field equations for the heat-conducting inviscid compressible fluid are obtained if balance laws~\eqref{eq:fieldeqn} are combined with the simplified closure relation for the stress tensor
\begin{alignat}{4}
	\mathrm{T}&=-\mathrm{p}\mathrm{I}\label{cr:ptensor}.
\end{alignat}
The closure relation~\eqref{cr:ptensor} is obtained from~\eqref{cr:CauchyStr} by setting the bulk and dynamic viscosity coefficient to zero, $\eta=\zeta=0$, from which it follows that the viscosity part~\eqref{cr:viscostress} of the stress tensor $\mathrm{T}$ \eqref{cr:CauchyStr} vanishes identically, $\upsigma\equiv 0$. The closure relation for $q$ is given by \eqref{eq:Fourier}. The resulting field equations are
\begin{subequations}\label{eq:fieldSimpl}
	\begin{align}	
	%%%%%% evolution equation for mass density
	\label{simplebal:mass}
	&\partial_t\rho
	+\diver\left(\rho{v}\right)=0,\\
	%%%%%% evolution equation for linear momentum density
	\label{simplebal:mom}
	&\partial_t(\rho{v})
	+\diver\left(\rho{v}\otimes{v}\right) =
	- \nabla \mathrm{p},\\
	%%%%%% evolution equation for internal energy density
	\label{simplebal:inten}
	&\partial_t(\rho\epsilon) + \diver\left(\rho\epsilon{v}\right) =
	-\diver q
	- \mathrm{p} \diver v.
	\end{align}
\end{subequations}
We show that the system of operator equations~\eqref{sys:pH1} encodes the system~\eqref{eq:fieldSimpl} in a weak form, by taking the example of field equation~\eqref{simplebal:inten}. For this we define the functional
\begin{alignat}{4}\label{eq:testfunctional}
&F(\zhet):=\int_{\Om}\phi(x)u(x)\dx &\quad&\text{with}&\quad& \phi\in C^\infty(\bar{\mathrm{\Omega}}).
\end{alignat}
Then the total functional derivative of \eqref{eq:testfunctional} is given by
\begin{alignat}{6}\label{testfunc:2}
&\vardif{F}{\zhet} = \begin{bmatrix}
\vardif{F}{\rho} & \vardif{F}{{M}}  & \vardif{F}{u}
\end{bmatrix}^\top\!\!=
\,\,\Big[\begin{array}{lcr}\!
0~~& 0 &~ \phi\!\!
\end{array}\Big]^\top.
\end{alignat}
The time evolution of the functional $F$ is described by
\begin{align}\label{eq:dtF}
\!\!\!\!\!\frac{\mathrm{d}F}{\mathrm{d}t}&\!=\!\aweak{\vardif{F}{\zhet}}{\dot{\zhet}}
\!\overset{\eqref{eq:ExOp1}}{=}\!\aweak{\vardif{F}{\zhet}}{\!\!\Big[\mathcal{J}\!-\mathcal{CDC}^\ast\Big]\widevardif{\mathcal{E}}{\zhet} + \mathcal{B}\mathbf{u}}\notag\\
&\!=\!\aweak{\vardif{F}{\zhet}}{\!\mathcal{J}\widevardif{\mathcal{E}}{\zhet}}\!-\!\aweak{\vardif{F}{\zhet}}{\!\mathcal{C}\mathcal{D}\mathcal{C}^\ast\widevardif{\mathcal{E}}{\zhet}}+\aweak{\vardif{F}{\zhet}}{\!\mathcal{B}\mathbf{u}}\!.
\end{align}
We consider the terms on the right-hand side of~\eqref{eq:dtF}. Under the assumption of a smooth solution and integration by parts (i.b.p.) the first term is rewritten such that
\begin{align}\label{id:1}
&\aweak{\vardif{F}{\zhet}}{\!\mathcal{J}\widevardif{\mathcal{E}}{\zhet}}\overset{\eqref{eq:JHydrodyn}}{=}\aweak{\vardif{F}{u}}{\!\op{J}{}{ u , {M}}\widevardif{\mathcal{E}}{M}}\notag\\
&\overset{\eqref{eq:JuM}}{=}
\int_\Om u \left(\widevardif{\mathcal{E}}{M} \cdot\nabla\right)\! \vardif{F}{u} +
\left(\widevardif{\mathcal{E}}{M}\cdot\nabla\right)\!\left(\vardif{F}{u}\mathrm{p}\right) \dx\notag\\
&\begin{aligned}[b]\hspace*{-1em}\overset{\text{i.b.p.}\,\&\, \eqref{eq:Epvardifs}}{=}&-\int_\Om  \vardif{F}{u} \Big[\diver\left(uv\right)+\mathrm{p}\diver v\Big]\dx\\
&+\int_{\partial\Om} \vardif{F}{u}\Big[u+\mathrm{p}\Big]v\cdot\nu\,\dS,\end{aligned}
\end{align}
where we have used~\eqref{eq:JuM} followed by partial integration. Also, we have used the relation $\wideminivardif{\mathcal{E}}{M}=(M/\rho)=v$, contained in the block vector~\eqref{eq:Epvardifs}. For the second term we obtain
\begin{align}\label{id:2}
&-\aweak{\vardif{F}{\zhet}}{\!\mathcal{C}\mathcal{D}\mathcal{C}^\ast\widevardif{\mathcal{E}}{\zhet}}
\overset{\eqref{eq:Fflux}}{=}\aweak{\!\vardif{F}{\zhet}}{\!\mathcal{C}F_\text{lux}}
=\aweak{\mathcal{C}^\ast\!\vardif{F}{\zhet}}{\!F_\text{lux}}\notag\\
&\overset{\eqref{eq:Cadj}}{=}
\int_\Om \nabla\!\left(\vardif{F}{u}\right)\cdot F_\text{lux}\, \dx
\overset{\text{i.b.p.}}{=}-\int_\Om \vardif{F}{u} \diver(F_\text{lux}) \, \dx \notag\\
&\hspace*{1cm}+\int_{\partial\Om}\vardif{F}{u}F_\text{lux}\cdot\nu\,\dS,
\end{align}
where $F_\text{lux}$ is defined in Equation~\eqref{eq:Fflux}. In our simplified setting, this is identical with the heat-flux vector $q$. The boundary contributions are reflected in the last term on the right-hand side of \eqref{eq:dtF}, given by
\begin{align}\label{id:3}
&\aweak{\vardif{F}{\zhet}}{\!\mathcal{B}\mathbf{u}}
\overset{\eqref{Bu:exergy}}{=}-\int_{\partial \Om}\vardif{F}{u}\left(\Big[u+\mathrm{p}\Big]v+F_\text{lux}\right)\cdot\nu\,\dS.
\end{align}
Expressing the terms on the right-hand side of time evolution equation~\eqref{eq:dtF} through~\eqref{id:1},~\eqref{id:2} and~\eqref{id:3}, results in
\begin{align}\label{eq:timeF}
\!\!\!\!\frac{\mathrm{d}F}{\mathrm{d}t}=-\int_\Om \vardif{F}{u} \Big[\diver\left(uv\right)+\mathrm{p}\diver v+ \diver(F_\text{lux})\Big]\dx.
\end{align}
On the other hand, the time evolution of the functional $F$ can also be calculated via
\begin{align}\label{eq:timeEvF2}
\frac{\mathrm{d}F}{\mathrm{d}t}=\int_{\Om}\vardif{F}{\zhet}\cdot\tpd{\zhet}{t}\dx\overset{\eqref{testfunc:2}}{=}\int_{\Om}\vardif{F}{u}\tpd{u}{t}\dx.
\end{align}
Combining~\eqref{eq:timeF} with \eqref{eq:timeEvF2} and using \eqref{testfunc:2}, we obtain
\begin{align}
\!\!\!\!\!0&=\!\!\int_\Om \vardif{F}{u} \Big[\diver\left(uv\right)+\mathrm{p}\diver v + \diver(F_\text{lux})+\partial\indices{_t}u\Big]\dx\notag\\
&\!\!\overset{\eqref{testfunc:2}}{=}\!\!\!\int_\Om\!\phi\Big[\diver\left(uv\right)+\mathrm{p}\diver v+ \diver(F_\text{lux})+\partial\indices{_t}u\Big]\dx,\!\!\!\!\!
\end{align}
and since $\phi\in C^\infty(\bar{\mathrm{\Omega}})$ is arbitrary, this gives
\begin{align}\label{eq:result}
\partial\indices{_t}u + \diver\left(u v\right)=-\diver\left(F_\text{lux}\right)-\mathrm{p}\diver v,
\end{align}
where $F_\text{lux}$ is the heat flux vector $q$ defined by \eqref{eq:Fflux}.
\section{Conclusion}
We have been able to rewrite the state space model of classical hydrodynamics from its representation in the GENERIC framework for open systems into a representation that may be seen as representing a generalized port-Hamiltonian system. The resulting single generator infinite dimensional state space representation is modeled with the same operators that are used in the original GENERIC formulation of classical hydrodynamics. This was achieved through the introduction of an exergy-like energy functional as a generating potential which allowed us to chose the internal energy density as one of the independent state variables. The skew-adjoint structure operator and self-adjoint dissipation operator of the resulting single generator state space representation satisfy crucial degeneracy requirements with respect to the total functional derivative of the exergy-like energy functional, which are inherited from the non-interacting conditions in the original GENERIC formulation.

Through a factorization of the dissipation operator we have been able to set up an extended skew-adjoint operator which is conjectured to define a Dirac structure. This is a focus of ongoing work.

\bibliography{References_11Nov18}
\bibliographystyle{apa}

\end{document}